\documentclass[a4paper]{amsart}
\usepackage{amssymb}
\usepackage{ifthen}
\usepackage{graphicx}

\newtheorem{thm}{Theorem}
\newtheorem{cor}{Corollary}
\newtheorem{lem}{Lemma}

\newtheorem{rem}{Remark}

\theoremstyle{definition}

\newcommand{\A}{{\mathcal A}}

\newcommand{\es}{{\mathcal S}}
\newcommand{\K}{{\mathcal K}}
\newcommand{\R}{{\mathcal R}}
\newcommand{\D}{{\mathbb D}}
\newcommand{\real}{{\operatorname{Re\,}}}

\begin{document}
\bibliographystyle{amsplain}

\title[differential inequalities and criteria for univalency]{Some differential inequalities and criteria for univalency in the unit disk}

\author[M. Obradovi\'{c}]{Milutin Obradovi\'{c}}
\address{Department of Mathematics,
Faculty of Civil Engineering, University of Belgrade,
Bulevar Kralja Aleksandra 73, 11000, Belgrade, Serbia.}
\email{obrad@grf.bg.ac.rs}

\author[N. Tuneski]{Nikola Tuneski}
\address{Department of Mathematics and Informatics, Faculty of Mechanical Engineering, Ss. Cyril and
Methodius
University in Skopje, Karpo\v{s} II b.b., 1000 Skopje, Republic of North Macedonia.}
\email{nikola.tuneski@mf.edu.mk}

\subjclass[2020]{30C45}
\keywords{differential inequalities, univalent functions, cirteria}

\begin{abstract}
In this paper, Jack lemma is used for obtaining several differential inequalities over analytic functions that later on, lead to new criteria for univalency in the unit disk.
\end{abstract}

\maketitle

\section{Introduction and preliminaries}

Let $\A$ denote the class of functions $f$ analytic in the open unit disk $\D=\{z:|z|<1\}$ and normalized such that $f(0)=f'(0)-1=0$, i.e.,
\begin{equation}\label{e1}
  f(z)=z+a_2z^2+a_3z^3+\cdots.
\end{equation}
Class $\es\subset \A$ consists of univalent (one-on-one and onto) functions, and its subclasses of starlike and convex functions of order $\alpha$, $0\le\alpha<1$, are defined respectively by
\begin{equation}\label{eq-2}
\begin{split}
\es^*(\alpha) &= \left\{ f\in\A : \real\left[\frac{zf'(z)}{f(z)}\right]>\alpha,\, z\in\D\right\},\\[2mm]
\K(\alpha) &= \left\{ f\in\A : \real\left[1+\frac{zf''(z)}{f'(z)}\right]>\alpha,\, z\in\D\right\}.
\end{split}
\end{equation}
For $\alpha=0$ we receive the classes of starlike and convex functions, $\es^*\equiv \es^*(0)$ and $\K\equiv\K(0)$. The class of starlike functions is very wide containing most of the other classes of univalent functions (convex, for example), but not all. One class that is not contained, nor contain, the class $\es^*$ is the class of functions with bounded turning
\begin{equation}\label{eq-3}
\R = \left\{ f\in\A : \real f'(z)>0,\, z\in\D\right\}.
\end{equation}
The geometrical characterisations of these classes correspond to their names. More details can be found in \cite{duren,book}.

\medskip

In this paper we will obtain several differential inequalities over a function $p$ analytic in the unit disk and such that $p(0)=1$, implying that it is a Caratheodory function, i.e., that $\real p(z)>0$ for all $z\in\D$. Afterwards, varying the function $p$ we will receive new criteria for a function $f\in\A$ to be in some of the classes of univalent functions defined above. This is a commonly used strategy for obtaining criteria for univalency of a function normalized and analytic in the unit disk. Collection of the most important results in that direction is given in \cite{mm}.

\medskip

For obtaining the results over the differential inequalities we will make use of the well known Jack lemma (\cite{jack}).

\medskip

\begin{lem}
Let $\omega$ be a nonconstant and analytic function in $\D$ with $\omega(0)=0$. If $|\omega|$ attains its maximum value on the circle $|z|=r$ at $z_0$, we have $z_0\omega'(z_0)=k\omega(z_0)$, $k\ge1$.
\end{lem}

\medskip

\section{Main results and consequences}

\medskip

\begin{thm}\label{th-1}
Let $p$ be analytic function on the unit dick  $\D$, normalized such that $p(0)=1$. If
\[ \left| zp'(z)+p(z)+p^2(z)-2\right| <\frac52 \qquad (z\in\D), \]
then $\real p(z)>0$ for all $z\in\D$.
\end{thm}

\medskip

\begin{proof}
Let consider function $\omega$, such that
\begin{equation}\label{eq-4}
p(z) = \frac{1+\omega(z)}{1-\omega(z)}.
\end{equation}
It is analytic in the unit disk and $\omega(0)=0$. It is enough to show that $|\omega(z)|<1$ for all $z\in\D$.

\medskip

On a contrary, let assume that $|\omega(z_0)|=1$ for some $z_0\in\D$. Then, according to the Jack lemma $z\omega'(z_0)=k\omega(z_0)$ for some $k\ge1$, and further, with $\omega(z_0)=e^{i\theta}$ and $t=\cos\theta$,
\[
\begin{split}
\left| z_0p'(z_0)+p(z_0)+p^2(z_0)-2\right|
&= \left| \frac{2z_0\omega'(z_0)}{[1-\omega(z_0)]^2} + \frac{1+\omega(z_0)}{1-\omega(z_0)} + \left[\frac{1+\omega(z_0)}{1-\omega(z_0)}\right]^2 -2\right| \\
&= 2 \cdot \left| \frac{\omega(z_0)\left[k+3-\omega(z_0)\right]}{[1-\omega(z_0)]^2}  \right| = 2 \cdot \frac{\left|k+3-\omega(z_0)\right|}{|1-\omega(z_0)|^2} \\
&= 2\cdot \frac{\sqrt{(k+3)^2-2(k+3)\cos\theta+1}}{2-2\cos\theta} \\
&= \frac{\sqrt{(k+3)^2-2(k+3)t+1}}{1-t} \equiv \varphi(t,k).
\end{split}
\]
Since $k\ge1$ and $-1\le t\le1$, we have that $\varphi_k(t,k) = \frac{3+k-t}{(1-t)\sqrt{(k+3)^2-2(k+3)t+1}} >0, $ and so
\[ \left| z_0p'(z_0)+p(z_0)+p^2(z_0)-2\right| \ge \varphi(t,1) = \frac{\sqrt{17-8t}}{1-t}, \]
which is an increasing function of $t$ on the interval $[-1,1]$. Thus
\[ \left| z_0p'(z_0)+p(z_0)+p^2(z_0)-2\right| \ge \varphi(-1,1) = \frac52. \]

\medskip

This is in contradiction with the condition of the theorem, so $|\omega(z)|<1$ for all $z\in\D$, implying $\real p(z)>0$ for all $z\in\D$.
\end{proof}

\medskip

For a function $f$ from $\A$, functions $\frac{zf'(z)}{f(z)}$, $1+\frac{zf''(z)}{f'(z)}$, $f'(z)$, and $\frac{f(z)}{z}$, are analytic in the unit disk and map the origin into 1. So, they can be chosen to be $p(z)$ in Theorem \ref{th-1}, and lead to the following results.

\medskip

\begin{cor}\label{cor-1}
Let $f\in\A$.
\begin{itemize}
  \item[($i$)] If
  \[\left| \frac{zf'(z)}{f(z)} \left[ 2 + \frac{zf''(z)}{f'(z)} \right] -2 \right| < \frac52 \qquad(z\in\D),\]
  then $f\in\es^*$.
  \item[($ii$)] If
  \[\left| \frac{z^2f'''(z)}{f'(z)} + 4\frac{zf''(z)}{f'(z)}  \right| < \frac52 \qquad(z\in\D),\]
  then $f\in\K$.
  \item[($iii$)] If
  \[\left| zf''(z) + f'(z) + f'^2(z) -2 \right| < \frac52 \qquad(z\in\D),\]
  then $\real f'(z)>0$ for all $z\in\D$.
  \item[($iv$)] If
  \[\left| f'(z) + \left[ \frac{f(z)}{z} \right]^2 -2 \right| < \frac52 \qquad(z\in\D),\]
  then $\real\frac{f(z)}{z}>0$ for all $z\in\D$.
\end{itemize}
\end{cor}

\medskip

In a similar way as Theorem \ref{th-1} we prove the following result.

\medskip

\begin{thm}\label{th-2}
Let $p$ be analytic function on the unit dick  $\D$, normalized such that $p(0)=1$. If
  \[ \left| zp'(z)+p(z)-p^2(z)\right| <\frac12 \qquad (z\in\D), \]
  then $\real p(z)>0$ for all $z\in\D$.
\end{thm}

\medskip

\begin{proof}
With the same assumptions and notations as in the proof of Theorem \ref{th-1}, and using similar technique, we again receive the conclusion by contradiction with the condition of the theorem:
\[
\begin{split}
\left| z_0p'(z_0)+p(z_0)-p^2(z_0)\right|
&= \left| \frac{2z_0\omega'(z_0)}{[1-\omega(z_0)]^2} + \frac{1+\omega(z_0)}{1-\omega(z_0)} - \left[\frac{1+\omega(z_0)}{1-\omega(z_0)}\right]^2 \right| \\
&= 2 \cdot \left| \frac{\omega(z_0)\left[k-1-\omega(z_0)\right]}{[1-\omega(z_0)]^2}  \right| \\
&= 2 \cdot \frac{\left|k-1-\omega(z_0)\right|}{|1-\omega(z_0)|^2}.
\end{split}
\]
The last is greater or equal to $\frac12$, with equality for $k=1$ and $\omega(z_0)=-1$.
\end{proof}

\medskip

As previously, by specifying functions $p$ in Theorem \ref{th-2} gives

\medskip

\begin{cor}\label{cor-2}
Let $f\in\A$.
\begin{itemize}
  \item[($i$)] If
  \[\left| \frac{zf'(z)}{f(z)} \left[ 2-2\frac{zf'(z)}{f(z)}+\frac{zf''(z)}{f'(z)} \right] \right| < \frac12 \qquad(z\in\D),\]
  then $f\in\es^*$.
  \item[($ii$)] If
  \[\left| \frac{z^2f'''(z)}{f'(z)} - 2\left[ \frac{zf''(z)}{f'(z)} \right]^2 \right| < \frac12 \qquad(z\in\D),\]
  then $f\in\K$.
  \item[($iii$)] If
  \[\left| zf''(z) + f'(z) - f'^2(z) \right| < \frac12 \qquad(z\in\D),\]
  then $\real f'(z)>0$ for all $z\in\D$.
  \item[($iv$)] If
  \[\left| f'(z) - \left[ \frac{f(z)}{z} \right]^2 \right| < \frac12 \qquad(z\in\D),\]
  then $\real\frac{f(z)}{z}>0$ for all $z\in\D$.
\end{itemize}
\end{cor}

\medskip

It is easy to verify that for $p(z)=\frac{zf'(z)}{f(z)} $,
\[ zp'(z)+p(z)-p^2(z) = -zf(z)\left[ \frac{z}{f(z)}\right]'', \]
and so, Corollary \ref{cor-2}(i) can be rewritten in the following form.

\medskip

\begin{thm}
If $f\in\A$ and
\[\left| f(z)\left[ \frac{z}{f(z)}\right]'' \right| \le \frac12 \qquad(z\in\D),\]
then $f\in\es^*$.
\end{thm}

\medskip

\begin{rem}
It is interesting that the condition
\[ \left| f(z)\left[ \frac{z}{f(z)} \right]' \right|<1 \quad(z\in\D), \]
which is equivalent to $ \left| \frac{zf'(z)}{f(z)}-1  \right|<1, $ $z\in\D$, also implies $f\in\es^*$.
\end{rem}

\medskip

If we repeat the proof of Theorem \ref{th-2}, but instead of \eqref{eq-4} work with function $\omega$ defined by
\[p(z) = (1-\alpha)\frac{1+\omega(z)}{1-\omega(z)}+\alpha,\]
we receive

\medskip

\begin{thm}\label{th-3}
Let $p$ be analytic function on the unit dick  $\D$, normalized such that $p(0)=1$, and let $0\le\alpha<1$. If
  \[ \left| zp'(z)+ \frac{1+\alpha}{1-\alpha}p(z)- \frac{1}{1-\alpha}p^2(z) - \frac{\alpha}{1-\alpha}\right| < \frac{1-\alpha}{2} \qquad (z\in\D), \]
  then $\real p(z)>\alpha$ for all $z\in\D$.
\end{thm}

\medskip

\begin{rem}
From the previous theorem we can obtain appropriate criteria as those in Corollary \ref{cor-1} and Corollary \ref{cor-2}. For example, for $\alpha=\frac12$, we have that
  \[ \left| zp'(z)+3p(z)-2p^2(z)-1\right| <\frac14 \,\, (z\in\D) \quad \Leftrightarrow \quad \real p(z)>\frac12 \,\, (z\in\D), \]
  and for $p(z)=\frac{zf'(z)}{f(z)} $, we receive that
    \[\left| \frac{zf'(z)}{f(z)} \left[ 4 - 3\frac{zf'(z)}{f(z)}+\frac{zf''(z)}{f'(z)} \right] - 1\right| < \frac14 \qquad(z\in\D),\]
  implies $f\in\es^*(1/2)$.
\end{rem}

\end{document}